\documentclass[11pt]{article}
\usepackage{latexsym}
\newcommand{\eproof}{\mbox{\ }\hfill $\Box$ \par \vskip 10pt}

\newtheorem{Theorem}{Theorem}[section]
\newtheorem{lemma}[Theorem]{Lemma}
\newtheorem{prop}[Theorem]{Proposition}

\newtheorem{corol}[Theorem]{Corollary}

\begin{document}

\title{Low frequency dispersive estimates for the wave equation
 in higher dimensions}

\author{{\sc Simon Moulin}}

\date{}

\maketitle

\abstract{We prove dispersive estimates at low frequency in
dimensions  $n\ge 4$ for the wave equation 
for a very large class of real-valued potentials, provided the zero
is neither an eigenvalue nor a resonance. This class includes
potentials $V\in L^\infty({\bf R}^n)$ satisfying
$V(x)=O\left(\langle x\rangle^{-(n+1)/2-\epsilon}\right)$,
$\epsilon>0$.}

\setcounter{section}{0}
\section{Introduction and statement of results}

High frequency dispersive estimates with loss of $(n-3)/2$ have been
recently proved in \cite{kn:V} for the wave equation with a
real-valued potential $V\in L^\infty({\bf R}^n)$, $n\ge 4$,
satisfying
$$|V(x)|\le C\langle x\rangle^{-\delta},\quad\forall x\in {\bf
R}^n,\eqno{(1.1)}$$ with constants $C>0$, $\delta>(n+1)/2$. The
problem of proving dispersive estimates at low frequency, however,
left open. The purposes of the present paper is to address
this problem. Such low frequency dispersive estimates for the
Schr\"odinger group have been recently proved in \cite{kn:MV} for a
large class of real-valued potentials (not necessarily in
$L^\infty$), and in particular for potentials satisfying (1.1) with
$\delta>(n+2)/2$. 

Denote by $G_0$ and $G$ the self-adjoint realizations of the
operators $-\Delta$ and $-\Delta+V$ on $L^2({\bf R}^n)$,
respectively. It is well known that, under the condition (1.1), the
absolutely continuous spectrums of the operators $G_0$ and $G$
coincide with the interval $[0,+\infty)$, and that $G$ has no
embedded strictly positive eigenvalues nor strictly positive
resonances. However, $G$ may have in general a finite number of
non-positive eigenvalues and that the zero may be a resonance. We
will say that the zero is a regular point for $G$ if it is neither
an eigenvalue nor a resonance in the sense that the operator
$1-V\Delta^{-1}$ is invertible on $L^1$ with a bounded inverse denoted by $T$. 
Let $P_{ac}$ denote the spectral projection onto the
absolutely continuous spectrum of $G$. Given any $a>0$, set
$\chi_a(\sigma)=\chi_1(\sigma/a)$, where $\chi_1\in C^\infty({\bf
R})$, $\chi_1(\sigma)=0$ for $\sigma\le 1$, $\chi_1(\sigma)=1$ for
$\sigma\ge 2$. Set $\eta_a=\chi(1-\chi_a)$, where $\chi$ denotes the
characteristic function of the interval $[0,+\infty)$. Clearly,
$\eta_a(G)+\chi_a(G)=P_{ac}$. As in the case of
the Schr\"odinger group (see \cite{kn:MV}), the dispersive estimates
for the low frequency part $e^{it\sqrt{G}}\eta_a(G)$, $a>0$ small,
turn out to be easier to prove when $n\ge 4$, and this can be done
for a larger class of potentials. In the present paper we will do so
for potentials satisfying
$$\sup_{y\in {\bf R}^n}\int_{{\bf
R}^n}\left(|x-y|^{-n+2}+|x-y|^{-(n-1)/2}\right)|V(x)|dx\le
C<+\infty.\eqno{(1.2)}$$ Clearly, (1.2) is fulfilled for potentials
satisfying (1.1). Our main result is the following

\begin{Theorem} Let $n\ge 4$, let $V$ satisfy (1.2) and assume that 
the zero is a regular point for $G$. Then, there exists a constant $a_0>0$ so that
for every $0<a\le a_0$, $0<\epsilon\ll 1$, $t$, we have the estimates
$$\left\|e^{it\sqrt{G}}G^{-(n+1)/4}\eta_a(G)\right\|_{L^1\to L^\infty}\le
C\langle t\rangle^{-(n-1)/2}\log(|t|+2),\eqno{(1.3)}$$
$$\left\|e^{it\sqrt{G}}G^{-(n+1)/4+\epsilon}\eta_a(G)\right\|_{L^1\to L^\infty}\le
C_\epsilon\langle t\rangle^{-(n-1)/2}.\eqno{(1.4)}$$
Moreover, for every $2\le p<+\infty$, we have the estimate
$$\left\|e^{it\sqrt{G}}G^{-\alpha(n+1)/4}\eta_a(G)\right\|_{L^{p'}\to L^p}\le
C\langle t\rangle^{-\alpha(n-1)/2},\eqno{(1.5)}$$ where
$1/p+1/p'=1$, $\alpha=1-2/p$, provided the operator $T$ is bounded on $L^{p'}$.
\end{Theorem}

\noindent
 {\bf Remark 1.} Note that our proof of the above estimates  
works out in the case $n=3$, too,
for potentials satysfying (1.2) as well as the condition 
$V\in L^{3/2-\epsilon}$ with some $0<\epsilon\ll 1$. In this case, however, 
a similar result has
been already proved by D'ancona and Pierfelice \cite{kn:DP}. In fact, in
\cite{kn:DP} the whole range of frequencies has been treated for a very large
subset of Kato potentials. 

Combining Theorem 1.1 with the estimates of \cite{kn:V}, we obtain the
following

\begin{corol}  Let $n\ge 4$, let $V$ satisfy (1.1) and assume that 
the zero is a
regular point for $G$. Then, for every  $2\le p<+\infty$,
$0<\epsilon\ll 1$, $t\neq 0$, we have the estimates
$$\left\|e^{it\sqrt{G}}G^{-(n+1)/4}\langle G\rangle^{-(n-3)/4-\epsilon}
P_{ac}\right\|_{L^1\to L^\infty}\le C_\epsilon
|t|^{-(n-1)/2}\log(|t|+2),\eqno{(1.6)}$$
$$\left\|e^{it\sqrt{G}}G^{-(n+1)/4+\epsilon}\langle G\rangle^{-(n-3)/4-2\epsilon}
P_{ac}\right\|_{L^1\to L^\infty}\le C_\epsilon
|t|^{-(n-1)/2},\eqno{(1.7)}$$
$$\left\|e^{it\sqrt{G}}G^{-\alpha(n+1)/4}\langle G\rangle^{-\alpha(n-3)/4}
P_{ac}\right\|_{L^{p'}\to L^p}\le
C|t|^{-\alpha(n-1)/2},\eqno{(1.8)}$$ where $1/p+1/p'=1$,
$\alpha=1-2/p$. Moreover, for every $0\le q\le (n-3)/2$, $2\le
p<\frac{2(n-1-2q)}{(n-3-2q)}$, we have
$$\left\|e^{it\sqrt{G}}G^{-\alpha(n+1)/4}\langle G\rangle^{-\alpha q/2}
P_{ac}\right\|_{L^{p'}\to L^p}\le
C|t|^{-\alpha(n-1)/2}.\eqno{(1.9)}$$
\end{corol}

Note that when $n=2$ and $n=3$ similar dispersive estimates (without loss
of derivatives) for the high frequency part $e^{it\sqrt{G}}\chi_a(G)$
are proved in \cite{kn:CCV} for potentials satisfying (1.1) (see also 
\cite{kn:GeV}, \cite{kn:DP}). 
For higher dimensions Beals \cite{kn:B} proved optimal
(without loss of derivatives) dispersive estimates for potentials
belonging to the Schwartz class. It seems that to avoid the loss of derivatives
in dimensions $n\ge 4$ one needs to impose some regularity condition on
the potential. Similar phenomenon also occurs in the case of the 
Schr\"odinger equation (see \cite{kn:GV}). Note that dispersive estimates 
without loss of derivatives for the Schr\"odinger group $e^{itG}$
in dimensions $n\ge 4$ are proved in \cite{kn:JSS} under the 
regularity condition $\widehat V\in L^1$. This result has been recently 
extended in \cite{kn:MV} to potentials $V$ satisfying (1.1) with $\delta>n-1$
as well as $\widehat V\in L^1$. 

To prove Theorem 1.1 we adapt the approach of \cite{kn:MV} to the wave equation.
It consists of proving uniform $L^1\to L^\infty$ dispersive estimates for
the operator  $e^{it\sqrt{G}}\psi(h^2G)$, where $\psi\in C_0^\infty
((0,+\infty))$, $h\gg 1$. To do so, we use Duhamel's formula for 
the wave equation (which in our case takes the form (2.12)). It turns out that
when $n\ge 4$ one can absorb the remaining terms taking the parameter $h$ 
big enough, so one does not need anymore to work on weighted $L^2$ spaces
(as in \cite{kn:V}). This allows to cover a larger class of potentials
not necessarily in $L^\infty$.

\section{Proof of Theorem 1.1}

Let $\psi\in C_0^\infty((0,+\infty))$. The following proposition is
proved in \cite{kn:MV} and that is why we omit the proof.

\begin{prop} Under the assumptions of Theorem 1.1,
there exist positive constants $C,\beta$ and $h_0$ so that the
following estimates hold
$$\left\|\psi(h^2G_0)\right\|_{L^1\to L^1}\le C,\quad
h>0,\eqno{(2.1)}$$
$$\left\|\psi(h^2G)\right\|_{L^1\to L^1}\le C,\quad
h\ge h_0,\eqno{(2.2)}$$
$$\left\|\psi(h^2G)-\psi(h^2G_0)T\right\|_{L^1\to L^1}\le Ch^{-\beta},\quad
h\ge h_0,\eqno{(2.3)}$$ where the operator
$$T=\left(1-V\Delta^{-1}\right)^{-1}:L^1\to L^1
\eqno{(2.4)}$$
is bounded by assumption.
\end{prop}

Set
$$\Phi(t,h)=e^{it\sqrt{G}}\psi(h^2G)-T^*e^{it\sqrt{G_0}}\psi(h^2G_0)T.$$
We will first show that Theorem 1.1 follows from the following

\begin{prop} Under the assumptions of Theorem 1.1, there exist
positive constants $C$, $h_0$ and $\beta$ so that for all $h\ge h_0$, $t$, we
have
$$\left\|\Phi(t,h)\right\|_{L^1\to L^\infty}\le
Ch^{-(n+1)/2-\beta}\langle t\rangle^{-(n-1)/2}.\eqno{(2.5)}$$
\end{prop}

By interpolation between (2.5) and the trivial bound
$$\left\|\Phi(t,h)\right\|_{L^2\to L^2}\le C,\eqno{(2.6)}$$
we obtain
$$\left\|\Phi(t,h)\right\|_{L^{p'}\to L^p}\le
Ch^{-\alpha(n+1)/2-\alpha\beta}\langle t\rangle^{-\alpha(n-1)/2},\eqno{(2.7)}$$ for every  $2\le p\le +\infty$, where $1/p+1/p'=1$,
$\alpha=1-2/p$. Now, writing
$$\sigma^{-\alpha(n+1)/4}\eta_a(\sigma)=\int_{a^{-1}}^\infty\psi(\sigma\theta)
\theta^{\alpha(n+1)/4}\frac{d\theta}{\theta}, \quad\sigma>0,$$
where $\psi(\sigma)=\sigma^{1-\alpha(n+1)/4}\chi'_1(\sigma)\in
C_0^\infty((0,+\infty))$, and using (2.7) we get (for $2<p\le
+\infty$)
$$\left\|e^{it\sqrt{G}}G^{-\alpha(n+1)/4}\eta_a(G)-T^*e^{it\sqrt{G_0}}
G_0^{-\alpha(n+1)/4}\eta_a(G_0)T \right\|_{L^{p'}\to L^p}$$ $$\le
\int_{a^{-1}}^\infty\left\|\Phi(t,\sqrt{\theta})\right\|_{L^{p'}\to
L^p}\theta^{-1+\alpha(n+1)/4}d\theta$$ $$\le
C\langle t\rangle^{-\alpha(n-1)/2}\int_{a^{-1}}^\infty\theta^{-1-\alpha\beta/2}d\theta
\le C\langle t\rangle^{-\alpha(n-1)/2},\eqno{(2.8)}$$ provided $a$ is taken small
enough. The estimate (1.5) follows from (2.8) and the fact that it holds
for $G_0$ (see \cite{kn:S}). Clearly, (1.3) follows from (2.8) with $p=+\infty$ and the estimate (A.1) in the appendix.
In the same way we get
$$\left\|e^{it\sqrt{G}}G^{-(n+1)/4+\epsilon}\eta_a(G)-T^*e^{it\sqrt{G_0}}
G_0^{-(n+1)/4+\epsilon}\eta_a(G_0)T \right\|_{L^1\to L^\infty}\le C\langle t\rangle^{-(n-1)/2},$$
which together with the estimate (A.2) in the appendix imply (1.4).\\

{\it Proof of Proposition 2.2.} We will derive (2.5) from the
following

\begin{prop} Under the assumptions of Theorem 1.1, there exist
positive constants $C$, $h_0$ and $\beta$ so that we have, for
 $\forall f\in L^1$,
$$\left\|e^{it\sqrt{G_0}}\psi(h^2G_0)f\right\|_{L^\infty}\le
Ch^{-(n+1)/2}\langle t\rangle^{-(n-1)/2}\|f\|_{L^1},\quad h\ge 1,\,\forall t,\eqno{(2.9)}$$
$$\int_{-\infty}^\infty
\left\|Ve^{it\sqrt{G_0}}\psi(h^2G_0)f
\right\|_{L^1}dt\le Ch^{-(n-1)/2}
\|f\|_{L^1},\quad h>0,\eqno{(2.10)}$$
$$\int_{-\infty}^\infty
\left\|Ve^{it\sqrt{G}}\psi(h^2G)f
\right\|_{L^1}dt\le Ch^{-1-\beta}
\|f\|_{L^1},\quad h\ge h_0.\eqno{(2.11)}$$
\end{prop}

We use Duhamel's formula
$$e^{it\sqrt{G}}=e^{it\sqrt{G_0}}+i\frac{\sin\left(t\sqrt{G_0}\right)}
{\sqrt{G_0}}
\left(\sqrt{G}-\sqrt{G_0}\right) -\int_0^t
\frac{\sin\left((t-\tau)\sqrt{G_0}\right)}{\sqrt{G_0}}Ve^{i\tau
\sqrt{G}}d\tau$$ to get the identity
$$\Phi(t;h)=\sum_{j=1}^2\Phi_j(t;h),\eqno{(2.12)}$$
where
$$\Phi_1(t;h)=\left(\psi_1(h^2G)-T^*\psi_1(h^2G_0)\right)e^{it\sqrt{G}}
\psi(h^2G)
$$ $$+T^*\psi_1(h^2G_0)e^{it\sqrt{G_0}}\left(\psi(h^2G)-\psi(h^2G_0)T\right)$$
 $$-iT^*\psi_1(h^2G_0)\sin\left(t\sqrt{G_0}\right)\left(\psi(h^2G)
-\psi(h^2G_0)T\right)
  $$ $$+iT^*\widetilde\psi_1(h^2G_0)\sin\left(t\sqrt{G_0}\right)\left(
 \widetilde\psi(h^2G)-\widetilde\psi(h^2G_0)T\right),$$
$$\Phi_2(t;h)=-h\int_0^tT^*\widetilde\psi_1(h^2G_0)\sin\left((t-\tau)
\sqrt{G_0}
\right)Ve^{i\tau\sqrt{ G}}\psi(h^2G)d\tau,$$ where $\psi_1\in
C_0^\infty((0,+\infty))$, $\psi_1=1$ on supp$\,\psi$,
$\widetilde\psi(\sigma)=\sigma^{1/2}\psi(\sigma)$,
$\widetilde\psi_1(\sigma)=\sigma^{-1/2}\psi_1(\sigma)$. Let $t>0$. By
Propositions 2.1 and 2.3, we have
$$\left\|\Phi_1(t;h)f\right\|_{L^\infty}\le
Ch^{-(n+1)/2-\beta}\langle t\rangle^{-(n-1)/2}\|f\|_{L^1}+Ch^{-\beta}
\left\|\Phi(t;h)f\right\|_{L^\infty},\eqno{(2.13)}$$
$$\langle t\rangle^{(n-1)/2}\left|\langle\Phi_2(t;h)f,g\rangle\right|$$ $$\le h
\int_0^{t/2}\langle t-\tau\rangle^{(n-1)/2} \left\|\sin\left((t-\tau)\sqrt{G_0}
\right)\widetilde\psi_1(h^2G_0)Tg\right\|_{
L^\infty}\left\|Ve^{i\tau \sqrt{G}}\psi(h^2G)f\right\|_{L^1}d\tau$$
 $$+h\int_{t/2}^t\left\|V\sin\left((t-\tau)\sqrt{G_0}
\right)\widetilde\psi_1(h^2G_0)Tg\right\|_{ L^1} \langle\tau\rangle^{(n-1)/2}
\left\|e^{i\tau \sqrt{G}}\psi(h^2G)f\right\|_{L^\infty}d\tau$$
 $$\le Ch^{-(n-1)/2}\|g\|_{ L^1}
\int_{-\infty}^{\infty}\left\|Ve^{i\tau\sqrt{
G}}\psi(h^2G)f\right\|_{L^1}d\tau$$
 $$+h\sup_{t/2\le\tau\le t}\langle\tau\rangle^{(n-1)/2} \left\|e^{i\tau\sqrt{
G}}\psi(h^2G)f\right\|_{L^\infty}
\int_{-\infty}^{\infty}\left\|V\sin\left((t-\tau)\sqrt{G_0}
\right)\widetilde\psi_1(h^2G_0)Tg\right\|_{L^1}d\tau$$
 $$\le Ch^{-(n+1)/2-\beta}\|g\|_{L^1}\|f\|_{L^1}
 +Ch^{-\beta}\|g\|_{L^1}\sup_{t/2\le\tau\le t}\langle\tau\rangle^{(n-1)/2} 
\left\|e^{i\tau\sqrt{
G}}\psi(h^2G)f\right\|_{L^\infty},$$ which clearly implies
 $$\langle t\rangle^{(n-1)/2}\left\|\Phi_2(t;h)f\right\|_{L^\infty}\le
Ch^{-(n+1)/2-\beta}\|f\|_{L^1}$$ $$+Ch^{-\beta}\sup_{t/2\le\tau\le
t}\langle\tau\rangle^{(n-1)/2} \left\|e^{i\tau\sqrt{
G}}\psi(h^2G)f\right\|_{L^\infty}.\eqno{(2.14)}$$ By (2.12)-(2.14),
we conclude
$$\langle t\rangle^{(n-1)/2}\left\|\Phi(t;h)f\right\|_{L^\infty}\le
Ch^{-(n+1)/2-\beta}\|f\|_{L^1}+Ch^{-\beta}
\langle t\rangle^{(n-1)/2}\left\|\Phi(t;h)f\right\|_{L^\infty}$$ $$
+Ch^{-\beta}\sup_{t/2\le\tau\le t}\langle\tau\rangle^{(n-1)/2}
\left\|\Phi(\tau;h)f\right\|_{L^\infty}.\eqno{(2.15)}$$ Taking $h$
big enough we can absorb the second and the third terms in the RHS
of (2.15), thus obtaining (2.5). Clearly, the case of $t<0$ can be
treated in the same way.\eproof

\section{ Proof of Proposition 2.3.} 

We will make use of the fact that
the kernel of the operator $e^{it\sqrt{ G_0}}\psi(h^2G_0)$ is of
the form $K_h(|x-y|,t)$, where
$$K_h(\sigma,t)=\frac{\sigma^{-2\nu}}{(2\pi)^{\nu+1}}\int_0^\infty
e^{it\lambda}{\cal J}_\nu(\sigma\lambda)\psi(h^2\lambda^2)\lambda
d\lambda= h^{-n}K_1(\sigma h^{-1},th^{-1}),\eqno{(3.1)}$$ where
${\cal J}_\nu(z)=z^\nu J_\nu(z)$,
$J_\nu(z)=\left(H_\nu^+(z)+H_\nu^-(z)\right)/2$ is the Bessel
function of order $\nu=(n-2)/2$. It is shown in \cite{kn:V}
(Section 2) that $K_h$ satisfies the estimates (for all $\sigma,
t>0$, $h\ge 1$)
$$\left|K_1(\sigma,t)\right|\le
C\langle t\rangle^{-s}\langle\sigma\rangle^{s-(n-1)/2},\quad\forall s\ge
0,\eqno{(3.2)}$$
$$\left|K_h(\sigma,t)\right|\le
Ch^{-(n+1)/2}\langle t\rangle^{-s}\sigma^{s-(n-1)/2},\quad 0\le s\le
(n-1)/2.\eqno{(3.3)}$$ Clearly, (2.9) follows from (3.3) with
$s=(n-1)/2$. It is not hard to see that (2.10) 
follows from (1.2) and the following

\begin{lemma} For all $\sigma, h>0$, $0\le s\le(n-1)/2$, we have
$$\int_{-\infty}^\infty |t|^s\left|K_h(\sigma,t)\right|dt\le
Ch^{-(n-1)/2}\sigma^{s-(n-1)/2}.\eqno{(3.4)}$$
\end{lemma}

{\it Proof.} In view of (3.1), it suffices to show (3.4) with
$h=1$. When $0<\sigma\le 1$, this follows from (3.2). Let now
$\sigma\ge 1$. We will use the fact that the function ${\cal J}_\nu$
can be decomposed as ${\cal
J}_\nu(z)=e^{iz}b_\nu^+(z)+e^{-iz}b_\nu^-(z)$, where $b_\nu^\pm(z)$
are symbols of order $(n-3)/2$ for $z\ge 1$. Then, we can decompose
the function $K_1$ as $K_1^++K_1^-$, where $K_1^\pm$ are defined by
replacing in the definition of $K_1$ the function ${\cal
J}_\nu(\sigma\lambda)$ by $e^{\pm
i\sigma\lambda}b_\nu^\pm(\sigma\lambda)$. Integrating by parts, we
get
$$\left|K_1^\pm(\sigma,t)\right|\le
C_m\sigma^{-(n-1)/2}|t\pm \sigma|^{-m},\eqno{(3.5)}$$ for every
integer $m\ge 0$. By (3.5),
$$\int_{-\infty}^\infty|t|^s\left|K_1^\pm(\sigma,t)\right|dt\le
\sigma^s\int_{-\infty}^\infty\left|K_1^\pm(\sigma,t)\right|dt
+\int_{-\infty}^\infty|t\pm\sigma|^s\left|K_1^\pm(\sigma,t)\right|dt$$
$$\le C_m\sigma^{s-(n-1)/2}\int_{-\infty}^\infty |t\pm
\sigma|^{-m}dt+C_m\sigma^{-(n-1)/2}\int_{-\infty}^\infty |t\pm
\sigma|^{-m+s}dt\le C\sigma^{s-(n-1)/2},\eqno{(3.6)}$$ which clearly implies
(3.4) in this case. \eproof

To prove (2.11) we will use the formula
$$e^{it\sqrt{G}}\psi(h^2G)=(i\pi h)^{-1}\int_0^\infty e^{it\lambda}
\varphi_h(\lambda)\left(R^+(\lambda)-R^-(\lambda)\right)d\lambda,
\eqno{(3.7)}$$
where $\varphi_h(\lambda)=\varphi_1(h\lambda)$, $\varphi_1(\lambda)=
\lambda\psi(\lambda^2)$, and $R^\pm(\lambda)=(G-\lambda^2\pm i0)^{-1}$
satisfy the identity
$$R^\pm(\lambda)\left(1+VR^\pm_0(\lambda)\right)
=R^\pm_0(\lambda).\eqno{(3.8)}$$
Here $R^\pm_0(\lambda)$ denote the outgoing and incoming free resolvents
with kernels given in terms of the Hankel functions, $H_\nu^\pm$, of
order $\nu=(n-2)/2$ by the formula
$$[R^\pm_0(\lambda)](x,y)=\pm i4^{-1}(2\pi)^{-\nu}|x-y|^{-n+2}\,
{\cal H}_\nu^\pm(\lambda|x-y|),$$
where ${\cal H}_\nu^\pm(z)=z^\nu H_\nu^\pm(z)$ satisfy
$$\left|\partial_z^j
{\cal H}_\nu^\pm(z)\right|\le C\langle z\rangle^{(n-3)/2},\quad
\forall z>0,\, j=0,1,$$
$$\left|{\cal H}_\nu^\pm(z)-{\cal H}_\nu^\pm(0)
\right|\le Cz^{1/2}\langle z\rangle^{(n-4)/2},\quad\forall z>0.$$
It follows easily from these bounds and (1.2) that
$$\left\|VR^\pm_0(\lambda)\right\|_{L^1\to L^1}\le C,\quad 0<\lambda\le 1,
\eqno{(3.9)}$$
$$\left\|VR^\pm_0(\lambda)-VR^\pm_0(0)
\right\|_{L^1\to L^1}\le C\lambda^{1/2},\quad 0<\lambda\le 1.
\eqno{(3.10)}$$
Since $1+VR^\pm_0(0)=1-V\Delta^{-1}$ is invertible on $L^1$ by assumption
with a bounded inverse denoted by $T$, it follows from (3.10) that there
exists a constant $\lambda_0>0$ so that the operator $1+VR^\pm_0(\lambda)$
is invertible on $L^1$ for $0<\lambda\le\lambda_0$. In view of (3.8), we
have
$$\sum_\pm \pm VR^\pm(\lambda)=-\sum_\pm \pm \left(1+VR_0^\pm(\lambda)
\right)^{-1}=-\sum_\pm \pm T\left(1+(VR_0^\pm(\lambda)-VR_0^\pm(0))T
\right)^{-1}$$ $$=\sum_\pm \pm T(VR_0^\pm(\lambda)-VR_0^\pm(0))T
\left(1+(VR_0^\pm(\lambda)-VR_0^\pm(0))T\right)^{-1}.\eqno{(3.11)}$$
By (3.7) and (3.11), 
$$Ve^{it\sqrt{G}}\psi(h^2G)=(i\pi h)^{-1}\sum_\pm \pm\int_{-\infty}^\infty TVP_h^\pm(t-\tau)
U_h^\pm(\tau)d\tau,\eqno{(3.12)}$$
where
$$P_h^\pm(t)=\int_0^\infty e^{it\lambda}
\widetilde\varphi_h(\lambda)\left(R^\pm_0(\lambda)-R^\pm_0(0)\right)d\lambda,$$
$$U_h^\pm(t)=\int_0^\infty e^{it\lambda}
\varphi_h(\lambda)T\left(1+(VR_0^\pm(\lambda)-VR_0^\pm(0))
T\right)^{-1}d\lambda,$$
where $\widetilde\varphi_h(\lambda)=\widetilde\varphi_1(h\lambda)$, 
$\widetilde\varphi_1\in C_0^\infty((0,+\infty))$ is such that
$\widetilde\varphi_1=1$ on supp$\,\varphi_1$. The kernel of the operator
$P_h^\pm(t)$ is of the form $A_h^\pm(|x-y|,t)$, where
$$A_h^\pm(\sigma,t)=\pm i4^{-1}(2\pi)^{-\nu}\sigma^{-n+2}\,
\int_0^\infty e^{it\lambda}
\widetilde\varphi_h(\lambda)\left({\cal H}_\nu^\pm(\sigma\lambda))-
{\cal H}_\nu^\pm(0)\right)d\lambda=h^{1-n}A_1^\pm(\sigma/h,t/h).\eqno{(3.13)}
$$

\begin{lemma} For all $\sigma>0$, $h\ge 1$, we have
$$\int_{-\infty}^\infty \left|A^\pm_h(\sigma,t)\right|dt\le
Ch^{-1/2}\left(\sigma^{-n+5/2}+\sigma^{-(n-1)/2}\right).\eqno{(3.14)}$$
\end{lemma}

{\it Proof.} In view of (3.13), it suffices to prove (3.14) with $h=1$.
Consider first the case $0<\sigma\le 1$. Using the inequality
$$\|\widehat f\|_{L^1}\le C\sum_{j=0}^1\sup_\lambda\langle\lambda\rangle
\left|\partial_\lambda^jf(\lambda)\right|,$$
we get
$$\sigma^{n-2}\,\int_{-\infty}^\infty \left|A^\pm_1(\sigma,t)\right|dt\le
C\sup_{\lambda\in{\rm supp}\,\widetilde\varphi_1}\left(
\left|{\cal H}_\nu^\pm(\sigma\lambda)-{\cal H}_\nu^\pm(0)\right|+\sigma
\left|\partial_\lambda{\cal H}_\nu^\pm(\sigma\lambda)\right|\right)\le
C\sigma^{1/2},$$
which is the desired bound. Let now $\sigma\ge 1$. 
We have
$$A^\pm_1(\sigma,t)=K^\pm_1(\sigma,t)+c^\pm
\sigma^{-n+2}\,\int_0^\infty e^{it\lambda}
\widetilde\varphi_1(\lambda)d\lambda,$$
where $c^\pm$ are constants and $K^\pm_1$ are as in the proof of Lemma 3.1.
Hence, in this case, (3.14) (with $h=1$) follows from (3.6) (with $s=0$).
\eproof

By (3.12), (3.14) and (1.2), we have
$$\int_{-\infty}^\infty
\left\|Ve^{it\sqrt{G}}\psi(h^2G)f
\right\|_{L^1}dt\le Ch^{-1}\sum_\pm\int_{-\infty}^\infty
\int_{-\infty}^\infty \left\|VP_h^\pm(t-\tau)U_h^\pm(\tau)f\right\|_{L^1}d\tau dt$$
 $$\le Ch^{-1}\sum_\pm\int_{-\infty}^\infty
\int_{-\infty}^\infty \int_{{\bf R}^n}\int_{{\bf R}^n}|V(x)|\left|A_h^\pm(|x-y|,t-
\tau)\right|\left|U_h^\pm(\tau)f(y)\right|dxdyd\tau dt$$
$$\le Ch^{-1}\sum_\pm
\int_{{\bf R}^n}\int_{{\bf R}^n}|V(x)|\left(\int_{-\infty}^\infty
\left|A_h^\pm(|x-y|,
\tau)\right|d\tau\right)\left(\int_{-\infty}^\infty 
\left|U_h^\pm(\tau)f(y)\right|d\tau\right)dxdy $$
$$\le Ch^{-3/2}\sum_\pm
\int_{{\bf R}^n}\int_{{\bf R}^n}|V(x)|\left(|x-y|^{-n+5/2}+|x-y|^{-(n-1)/2}
\right)\int_{-\infty}^\infty 
\left|U_h^\pm(\tau)f(y)\right|d\tau dxdy $$
 $$\le Ch^{-3/2}\sum_\pm
\int_{{\bf R}^n}\int_{-\infty}^\infty 
\left|U_h^\pm(\tau)f(y)\right|d\tau dy.\eqno{(3.15)} $$
Thus, (2.11) follows from (3.15) and the following

\begin{lemma} There exists a constant $h_0>0$ so that for $h\ge h_0$
we have
$$\int_{{\bf R}^n}\int_{-\infty}^\infty 
\left|U_h^\pm(t)f(x)\right|dt dx\le C\|f\|_{L^1}.\eqno{(3.16)} $$
\end{lemma}

{\it Proof.} Using the identity
$$T\left(1+(VR_0^\pm(\lambda)-VR_0^\pm(0))T\right)^{-1}$$ $$=
T-T(VR_0^\pm(\lambda)-VR_0^\pm(0))T
\left(1+(VR_0^\pm(\lambda)-VR_0^\pm(0))T\right)^{-1},$$
we obtain
$$U_h^\pm(t)=T\widehat\varphi_h(t)-\int_{-\infty}^\infty TVP_h^\pm(t-\tau)
U_h^\pm(\tau)d\tau.\eqno{(3.17)}$$
Since
$$\int_{-\infty}^\infty |\widehat\varphi_h(t)|dt=
h^{-1}\int_{-\infty}^\infty |\widehat\varphi_1(t/h)|dt
=\int_{-\infty}^\infty |\widehat\varphi_1(t)|dt,$$
as above, we have
$$\int_{{\bf R}^n}\int_{-\infty}^\infty 
\left|U_h^\pm(t)f(x)\right|dt dx\le C\|f\|_{L^1}\int_{-\infty}^\infty 
|\widehat\varphi_h(t)|dt$$ $$+C\int_{-\infty}^\infty
\int_{-\infty}^\infty \int_{{\bf R}^n}\int_{{\bf R}^n}|V(x)|\left|A_h^\pm(|x-y|,t-
\tau)\right|\left|U_h^\pm(\tau)f(y)\right|dxdyd\tau dt$$
$$\le  C\|f\|_{L^1}+Ch^{-1/2}
\int_{{\bf R}^n}\int_{-\infty}^\infty 
\left|U_h^\pm(\tau)f(y)\right|d\tau dy,$$
which implies (3.16) provided $h$ is taken big enough.
\eproof

\appendix
\section{Appendix}

The following low frequency dispersive estimates for the free wave group are more or less known, but we will give a proof for the
sake of completeness. We have the following

\begin{prop} Let $n\ge 3$. Then for every $0<\epsilon\ll 1$, $t$, we have the estimates
$$\left\|e^{it\sqrt{G_0}}G_0^{-(n+1)/4}\eta_a(G_0)\right\|_{L^1\to L^\infty}\le C\langle t\rangle^{-(n-1)/2}\log\left(
|t|+2\right),\eqno{(A.1)}$$
$$\left\|e^{it\sqrt{G_0}}G_0^{-(n+1)/4+\epsilon}\eta_a(G_0)\right\|_{L^1\to L^\infty}\le C_\epsilon\langle t\rangle^{-(n-1)/2}.\eqno{(A.2)}$$
\end{prop}

{\it Proof.} The kernel of the operator in the LHS of (A.1) is of the form $K(|x-y|,t)$, where
$$K(\sigma,t)=c_n\sigma^{-n+2}\int_0^\infty e^{it\lambda}\lambda^{1-(n+1)/2}\eta_a(\lambda^2){\cal J}_\nu(\sigma\lambda)d\lambda.$$
When $|t|\le 2$, using that ${\cal J}_\nu(z)=O(z^{n-2})$, $\forall z>0$, we have $|K(\sigma,t)|\le Const$, which implies
(A.1) in this case. In what follows we will suppose $|t|\ge 2$. Let $\phi\in C_0^\infty({\bf R})$, $\phi(\mu)=1$ for $|\mu|\le 1$,
$\phi(\mu)=0$ for $|\mu|\ge 2$. We write $K=K_1+K_2$, where
$$K_1(\sigma,t)=c_n\sigma^{-n+2}\int_0^\infty e^{it\lambda}\lambda^{1-(n+1)/2}\eta_a(\lambda^2)(\phi{\cal J}_\nu)(\sigma\lambda)d\lambda,$$
$$K_2(\sigma,t)=c_n\sigma^{-n+2}\int_0^\infty e^{it\lambda}\lambda^{1-(n+1)/2}\eta_a(\lambda^2)((1-\phi){\cal J}_\nu)(\sigma\lambda)d\lambda.$$
Since $((1-\phi){\cal J}_\nu)(z)=O(z^{(n-3)/2})$, $\forall z>0$, we have
$$|K_2(\sigma,t)|\le C\sigma^{-(n-1)/2}\int_{\sigma^{-1}}^{Const}\lambda^{-1}d\lambda\le C\sigma^{-(n-1)/2}\log\langle\sigma\rangle.\eqno{(A.3)}$$ It follows from (A.3) that for $|t|/2\le\sigma\le 2|t|$, we have
$$|K_2(\sigma,t)|\le C|t|^{-(n-1)/2}\log|t|.\eqno{(A.4)}$$
Let now $\sigma\not\in [|t|/2,2|t|]$. We write $K_2$ as $K_2^++K_2^-$, where 
$$K_2^\pm(\sigma,t)=c_n\sigma^{-n+2}\int_0^\infty e^{i(t\pm\sigma)\lambda}\lambda^{1-(n+1)/2}\eta_a(\lambda^2)((1-\phi)b^\pm_\nu)(\sigma\lambda)d\lambda,$$
with functions $b_\nu^\pm$ satisfying
$$|\partial_z^jb_\nu^\pm(z)|\le C_jz^{(n-3)/2-j},\quad\forall j\ge 0,\,z\ge 1.$$
Integrating by parts $m\ge 1$ times we get
$$|K_2^\pm(\sigma,t)|\le C\sigma^{-n+2}|t\pm\sigma|^{-m}\int_0^\infty \sum_{j=0}^m\sigma^{m-j}\left|\partial^j_\lambda(
\lambda^{1-(n+1)/2}\eta_a(\lambda^2))\right|\left|(\partial^{m-j}_\lambda(1-\phi)b^\pm_\nu)(\sigma\lambda)\right|d\lambda$$
$$\le C\sigma^{-n+2}|t\pm\sigma|^{-m}\int_{\sigma^{-1}}^{Const} \sum_{j=0}^m\sigma^{m-j}
\lambda^{1-(n+1)/2-j}(\sigma\lambda)^{(n-3)/2-(m-j)}d\lambda$$
$$\le C\sigma^{-(n-1)/2}|t\pm\sigma|^{-m}\int_{\sigma^{-1}}^{Const}\lambda^{-1-m}d\lambda
\le C\sigma^{m-(n-1)/2}|t\pm\sigma|^{-m}\int_1^\infty\mu^{-1-m}d\mu$$
$$\le C\sigma^{m-(n-1)/2}|t\pm\sigma|^{-m}\le C_m\sigma^{m-(n-1)/2}|t|^{-m},\eqno{(A.5)}$$
since $|t\pm\sigma|\ge |t|/2$ in this case,
for all integers $m\ge 1$, and hence for all real $m\ge 1$. Taking $m=(n-1)/2$ in (A.5) we get
$$|K_2(\sigma,t)|\le C|t|^{-(n-1)/2},\quad\mbox{ if}\quad \sigma\not\in [|t|/2,2|t|].\eqno{(A.6)}$$
To deal with $K_1$ we will use that $(\phi{\cal J}_\nu)(z)=z^{n-2}g(z)$ with a function $g\in C_0^\infty({\bf R})$.
We write
$$K_1(\sigma,t)=c_n\int_0^\infty e^{it\lambda}\lambda^{(n-3)/2}\eta_a(\lambda^2)g(\sigma\lambda)d\lambda.$$

\begin{lemma} For every $k\ge 1$, we have
$$\left|\int_0^\infty e^{it\lambda}\lambda^{k-1}\eta_a(\lambda^2)g(\sigma\lambda)d\lambda\right|\le C_k|t|^{-k},\eqno{(A.7)}$$
with a constant $C_k>0$ indpendent of $t$ and $\sigma$.
\end{lemma}

{\it Proof.} If $k\ge 1$ is an integer, we integrate by parts $k$ times to get
$$\left|\int_0^\infty e^{it\lambda}\lambda^{k-1}\eta_a(\lambda^2)g(\sigma\lambda)d\lambda\right|$$ $$\le |t|^{-k}
\int_0^\infty \left|\partial_\lambda^k(\lambda^{k-1}\eta_a(\lambda^2)g(\sigma\lambda))\right|d\lambda+|t|^{-k}\left|
\partial_\lambda^{k-1}(\lambda^{k-1}\eta_a(\lambda^2)g(\sigma\lambda))|_{\lambda=0}\right|$$
 $$\le |t|^{-k}
\int_0^\infty\sum_{j=1}^k \sigma^j\lambda^{j-1}|(\partial^j_\lambda g)(\sigma\lambda)|d\lambda
+|t|^{-k}\int_0^\infty|\eta'_a(\lambda^2)||g(\sigma\lambda)|d\lambda
+|t|^{-k}|g(0)|$$
 $$\le |t|^{-k}\sum_{j=1}^k 
\int_0^\infty(\sigma\lambda)^{j-1}(\partial^j_\lambda g)(\sigma\lambda)|d(\sigma\lambda)
+|t|^{-k}\int_0^\infty|\eta'_a(\lambda^2)|d\lambda
+|t|^{-k}|g(0)|\le C_k|t|^{-k}.$$
For all real $k\ge 1$, (A.7) follows easily by complex interpolation. 
\eproof

Applying (A.7) with $k=(n-1)/2$ we get
$$|K_1(\sigma,t)|\le C|t|^{-(n-1)/2}.\eqno{(A.8)}$$
Now (A.1) follows from (A.4), (A.6) and (A.8).

To prove (A.2) observe that the function
$$\widetilde K_2(\sigma,t)=c_n\sigma^{-n+2}\int_0^\infty e^{it\lambda}\lambda^{1+2\epsilon-(n+1)/2}\eta_a(\lambda^2)((1-\phi){\cal J}_\nu)(\sigma\lambda)d\lambda,$$
satisfies the bound
$$|\widetilde K_2(\sigma,t)|\le C\sigma^{-(n-1)/2}\int_{\sigma^{-1}}^{Const}\lambda^{-1+2\epsilon}d\lambda\le C_\epsilon\sigma^{-(n-1)/2}.\eqno{(A.9)}$$
Hence, for $|t|/2\le\sigma\le 2|t|$, we have
$$|\widetilde K_2(\sigma,t)|\le C|t|^{-(n-1)/2}.\eqno{(A.10)}$$
The rest of the proof is exactly as above.
\eproof

Universit\'e de Nantes,
 D\'epartement de Math\'ematiques, UMR 6629 du CNRS,
 2, rue de la Houssini\`ere, BP 92208, 44332 Nantes Cedex 03, France

e-mail: simon.moulin@math.univ-nantes.fr

\end{document}